\documentclass[11pt]{article}
\usepackage{}
\usepackage{mathrsfs}
\usepackage{amsfonts}

\usepackage{CJK}
\usepackage{amsfonts,amssymb,amsmath,mathrsfs}

\usepackage{color,latexsym,amsfonts}
 \newtheorem{theorem}{Theorem}[section]
 \newtheorem{lemma}[theorem]{Lemma}
 \newtheorem{corollary}[theorem]{Corollary}
 \newtheorem{proposition}[theorem]{Proposition}
 
 \newtheorem{Definition}[theorem]{Definition}
 \newtheorem{remark}[theorem]{Remark}
 \newtheorem{condition}[theorem]{Condition}

 \def\btheorem{\begin{theorem}\sl{}\def\etheorem{\end{theorem}}}
 \def\bcorollary{\begin{corollary}\sl{}\def\ecorollary{\end{corollary}}}
 
 \def\bproposition{\begin{proposition}\sl{}\def\eproposition{\end{proposition}}}
 \def\bremark{\begin{remark}\sl{}\def\eremark{\end{remark}}}

 \def\beqlb{\begin{eqnarray}}\def\eeqlb{\end{eqnarray}}
 \def\beqnn{\begin{eqnarray*}}\def\eeqnn{\end{eqnarray*}}

 \def\<{\langle}\def\>{\rangle}

 \def\ar{&\!\!}

 \def\eqref#1{{\rm(\ref{#1})}}

\begin{document}

\

\noindent{}

\bigskip\bigskip

\centerline{\Large\bf Harnack Inequalities and Applications for}

\smallskip

\centerline{\Large\bf
 Stochastic Differential Equations Driven by}

\smallskip

\centerline{\Large\bf
Fractional Brownian Motion}

\smallskip
\
\bigskip\bigskip

\centerline{Xi-Liang Fan}

\bigskip

\centerline{School of Mathematical Sciences, Beijing Normal
 University,}

\centerline{Beijing 100875, People's Republic of China}

\smallskip

\bigskip\bigskip

{\narrower{\narrower

\noindent{\bf Abstract.} In the paper, Harnack inequalities are established for stochastic differential equations driven by fractional Brownian motion with Hurst parameter $H<\frac{1}{2}$. As applications, strong Feller property, log-Harnack inequality and entropy-cost inequality are given.}
}

\bigskip
 \textit{Mathematics Subject Classifications (2000)}: Primary 60H15

\bigskip

\textit{Key words and phrases}: Harnack inequality, stochastic differential equation, fractional Brownian motion.


\section{Introduction}

\setcounter{equation}{0}
Under a curvature condition, Wang \cite{Wang1} introduced dimensional-free Harnack inequality for diffusions on Riemannian manifold.
This type of inequality has been studied extensively, see, for example, Aida and Kawabi \cite{Aida1,Aida2} for infinite dimensional diffusion processes; Wang \cite{Wang2} for stochastic generalized porous media equations; R\"{o}ckner and Wang \cite{Rockner} for generalizes Mehler semigroup; \cite{Abdelhadi} for stochastic functional differential equation; Ouyang \cite{Ouyang} for Ornstein-Uhnelbeck processes and multivalued stochastic evolution equations etc.

Harnack inequality has various applications, see, for instance, \cite{Bobkov,Rockner,Rockner2,Wang6,Wang7} for strong Feller property and contractivity properties; \cite{Aida1,Aida2} for short times behaviors of infinite dimensional diffusions; \cite{Bobkov,Gong} for heat kernel estimates and entropy-cost inequalities. \cite{Aida1,Kawabi,Rockner,Wang1} established Harnack inequalities
using the method of derivative formula. In order to obtain Harnack inequality on manifolds with unbounded below curvatures, \cite{Arnaudon} introduced the approach of coupling and Girsanov transformations. In the paper, we will use the above two methods to establish Harnack inequalities for stochastic differential equations driven by fractional Brownian motion.

Solutions of the stochastic differential equations driven by fractional Brownian motion have been studied intensively in recent years, for example see \cite{Lyons,Nualart2} using the pathwise approach; see \cite{Coutin} using the tools of rough path analysis introduced in \cite{Lyons}. We prove Harnack inequality for stochastic differential equations driven by fractional Brownian motion with Hurst parameter $H<\frac{1}{2}$. As applications of the Harnack inequality, the strong Feller property and the log-Harnack inequality are derived. We also get the entropy-cost inequality with respect to the Euclidian distance.

The paper is organized as follows. In section 2, we give some preliminaries on fractional Brownian motion. Section 3 prove the Harnack inequality by
using the approach of coupling and Girsanov transformations,  and present their applications.  In section 4, we are devoted to establish derivative formula and give the corresponding Harnack inequality.

Harnack inequality

\section{Preliminaries}

\setcounter{equation}{0}
Let $B^H=\{B_t^H, t\in[0,T]\}$ be a fractional Brownian motion with Hurst parameter $H\in(0,1)$ defined on the probability
space $\Omega,\mathscr{F},\mathbb{P}$, i.e., $B^H$ is a centered Gauss process with the covariance function
\beqnn
 R_H(t,s)=\mathbb{E}(B_t^HB_s^H)=\frac{1}{2}\left(t^{2H}+s^{2H}-|t-s|^{2H}\right).
\eeqnn
In particular, if $H=\frac{1}{2}, B$ is a Brownian motion.
It is well known that if $H\neq\frac{1}{2}, B^H$ does not have independent increments
and has $\alpha$-order H\"{o}lder continuous path for all $\alpha\in(0,H)$.

For each $t\in[0,T]$, we denote by $\mathcal {F}_t$ the $\sigma$-algebra generated by the random variables $\{B_s^H:s\in[0,T]\}$ and the
$\mathbb{P}$-null sets.

We denote by $\mathscr{E}$ the set of step functions on $[0,T]$. Let $\mathcal {H}$ be the Hilbert space defined as the closure of $\mathscr{E}$ with respect to the scalar product
\beqnn
 \langle I_{[0,t]},I_{[0,s]}\rangle R_H(t,s)=\mathbb{E}(B_t^HB_s^H)=\frac{1}{2}\left(t^{2H}+s^{2H}-|t-s|^{2H}\right).
\eeqnn
The mapping $I_[0,T]\mapsto B_t^H$ can be extended to an isometry between $\mathcal {H}$ and the Gauss space
$\mathcal {H}_1$ associated with $B^H$. Denote this isometry by $\phi\mapsto B^H(\phi)$. For more details, one can see \cite{Nualart3}.
On the other hand, from \cite{Decr}, we know the covariance kernel $R_H(t,s)$ can be written as
\beqnn
 R_H(t,s)=\int_0^{t\wedge s}K_H(t,r)K_H(s,r)dr,
\eeqnn
where $K_H$ is a square integrable kernel given by
\beqnn
 K_H(t,s)=\Gamma(H+\frac{1}{2})^{-1}(t-e)^{H-\frac{1}{2}}F(H-\frac{1}{2},\frac{1}{2}-H,H+\frac{1}{2},1-\frac{t}{s}),
\eeqnn
in which $F(\cdot,\cdot,\cdot,\cdot)$ is the Gauss hypergeometric function.\\
Define the linear operator $K_H^*\mathscr{E}\rightarrow L^2[0,T]$ as follows
\beqnn
 (K_H^*\phi)(s)=K_H(t,s)\phi(s)+\int_s^T(\phi(r)-\phi(s))\frac{\partial K_H}{\partial r}(r,s)dr.
\eeqnn
By \cite{Alos}, we know that, for all $\phi,\psi\in\mathscr{E}$, $ \langle K_H^*\phi,K_H^*\psi\rangle_L^2[0,T]=\langle \phi,\psi\rangle$ holds.
From B.L.T. theorem, $K_H^*$ can be extended to an isometry between $\mathcal {H}$ and $L^2[0,T]$. Therefore, according to \cite{Alos},
the process $\{W_t=B((K_H^*)^{-1}(I_{[0,t]})),t\in[0,T]\}$ is a Wiener process, and $B^H$ has the following integral representation
\beqnn
 B^H_t=\int_0^tK_H(t,s)dW_s.
\eeqnn
By \cite{Decr}, the operator $K_H:L^2[0,T]\rightarrow I_{0+}^{H+\frac{1}{2}}(L^2[0,T])$ associated with the square integrable kernel $K_H(\cdot,\cdot)$ is defined as follows
\beqnn
 (K_Hf)(t):=\int_0^tK_H(t,s)f(s)ds,\ f\in L^2[0,T].
\eeqnn
where $ I_{0+}^{H+\frac{1}{2}}$  is the $\alpha$-order left fractional Riemannian-Liouville integral operator on ${0,T}$, one can see \cite{Samko}.
It is an isomorphism and for each $f\in L^2[0,T]$,
\beqnn
 (K_Hf)(s)=I_{0+}^{2H}s^{\frac{1}{2}-H}I_{0+}^{\frac{1}{2}-H}s^{H-\frac{1}{2}}f,\ H\leq\frac{1}{2},\\
 (K_Hf)(s)=I_{0+}^{1}s^{H-\frac{1}{2}}I_{0+}^{H-\frac{1}{2}}s^{\frac{1}{2}-H}f,\ H\geq\frac{1}{2}.
\eeqnn
As a consequence, for every $h\in I_{0+}^{H+\frac{1}{2}}(L^2[0,T])$, the inverse operator $K_H^{-1}$ is of the following form
\beqnn
 (K_H^{-1}h)(s)=s^{H-\frac{1}{2}}D_{0+}^{H-\frac{1}{2}}s^{\frac{1}{2}-H}h',\ H>\frac{1}{2},\\
 (K_H^{-1}h)(s)=s^{\frac{1}{2}-H}D_{0+}^{\frac{1}{2}-H}s^{H-\frac{1}{2}}D_{0+}^{2H}h,\ H<\frac{1}{2},
\eeqnn
where $D_{0+}^{H-\frac{1}{2}}(D_{0+}^{\frac{1}{2}-H})$ is $H-\frac{1}{2}(\frac{1}{2}-H)$-order left-sided Riemannian-Liouville derivative, one also can see \cite{Samko}.\\
In particular, if $h$ is absolutely continuous, we have
\beqnn
 (K_H^{-1}h)(s)=s^{H-\frac{1}{2}}I_{0+}^{\frac{1}{2}-H}s^{\frac{1}{2}-H}h',\ H<\frac{1}{2}.
\eeqnn
In \cite{Nualart}, D.Nualart and Y.Ouknine discussed the following stochastic differential equations driven by fractional Brownian motion on $\mathbb{R}$,
\beqlb\label{2.1}
dX_t=b(t,X_t)dt+dB_t^H,\  X_0=x.
 \eeqlb
They proved the existence and uniqueness of a strong solution for \eqref{2.1} when $b(t,x)$ is a Borel function with linear growth in $x$ in case $H\leq\frac{1}{2}$.\\
The aim of the paper is to consider
the Harnack inequality for the equation \eqref{2.1} in case $H<\frac{1}{2}$.\\
We define $P_tf(x):=\mathbb{E}f(X_t^x), \ t\in[0,T], \ f\in \mathscr{B}_b(\mathbb{R})$, where
$X_t^x$ is the solution to the equation \eqref{2.1} and
$\mathscr{B}_b(\mathbb{R})$ denotes the set of all bounded
measurable functions on $\mathbb{R}$.

\section{Main results and proofs}

\setcounter{equation}{0}

Let us start with the following hypothesis (H1):
\begin{itemize}
 \item[(i)]
 $|b(t,x)-b(t,y)|\leq K|x-y|$,  $\forall x,y\in\mathbb{R}, t\in
[0,T]$, where $K>0$ is a constant;
 \par
 \item[(ii)]The mapping $t\mapsto b(t,0)$ is bounded on $[0,T]$.
\end{itemize}
It is clear that under (H1), the equation \eqref{2.1} has a unique solution. Furthermore, we can give the Harnack inequality for the equation \eqref{2.1} as follows.
\btheorem\label{t3.1}  If {\rm(H1)} holds, then for any nonnegative $f\in\mathscr{B}_b(\mathbb{R})$ and $t>0, x,y\in\mathbb{R}$,
\beqnn
(P_Tf(y))^p\leq P_Tf^p(x){\rm exp}[\frac{p}{p-1}C(T,K,H)|x-y|^2],
\eeqnn
where
$C(T,K,H)=\left(\frac{B(\frac{3}{2}-H,\frac{1}{2}-H)}{\Gamma(\frac{1}{2}-H)}\right)^2\frac{T^{2-2H}}{K^{-2}(1-e^{-2KT})^2(1-H)}$.
\etheorem

\emph{Proof.} The proof will be divided into three steps.\\
Step 1. Consider the following coupled stochastic differential equation
\beqlb\label{3.1}
\ar dY_t\ar=b(t,Y_t)dt+dB_t^H+u_tdt,\ \ Y_0=y,\\
\ar dX_t\ar=b(t,X_t)dt+dB_t^H, \ \ X_0=x,
\eeqlb
where the drift term $u_t$ of the equation  \eqref{3.1}  is of the following form
\beqnn \eta_t\cdot\frac{X_t-Y_t}{|X_t-Y_t|}I_{\{t<\tau\}},\eeqnn
$\tau$ is the coupling time of $X_t$ and $Y_t$ defined by
\beqnn \tau=\inf\{t\geq 0: X_t=Y_t\},\eeqnn and $\eta_t$ is a deterministic function on $[0,\infty)$ specified later such that the force $u_t$ can
make the two processes $X$ and $Y$ move together before time $T$.\\
It is obvious that the assumption {\rm(H1)} implies $|b(t,x)|\leq C(1+|x|)$, then, according to \cite[Theorem 8]{Nualart},
the equation \eqref{3.1} has a unique solution.\\
Note that $d(X_t-Y_t)=\left(b(t,X_t)-b(t,Y_t)\right)dt-u_tdt$, thus applying the Tanaka formula to $|X_t-Y_t|$, we have for $t<\tau$
\beqnn
d|X_t-Y_t|\ar=\ar {\rm sgn}(X_t-Y_t)d(X_t-Y_t)\cr
\ar=\ar {\rm sgn}(X_t-Y_t)\left(b(t,X_t)-b(t,Y_t)\right)dt-\eta_tdt.
\eeqnn
By ${\rm(H1)}$, for all $t<\tau$ we get
\beqnn
d|X_t-Y_t|\leq (K|X_t-Y_t|-\eta_t)dt.
\eeqnn
This implies that
\beqlb\label{3.3}
e^{-K(T\wedge\tau)}|X_T-Y_T|\leq |x-y|+\int_0^Te^{-Kt}\eta_tdt.
\eeqlb
Choosing \beqnn\eta_t=\frac{e^{-Kt}}{\int_0^Te^{-2Kt}dt}\cdot|x-y|,\ \ t\geq0.\eeqnn
We conclude that $\tau\leq T$ and $X_T=Y_T, a.s.$ Otherwise, if $\tau>T$, by \eqref{3.3} we get $X_T=Y_T$. But this contradicts with the assumption
that $\tau>T$.\\
Step 2.
Let $\tilde{B}_t^H=\int_0^tu_sds+B_t^H,\ \forall t\in[0,T].$ By simple calculus, we know that $\int_0^Tu_t^2dt<\infty$. Hence, $\int_0^\cdot u_rdr\in I_{0+}^{H+\frac{1}{2}}(L^2([0,T]))$. According to integral representation of fractional Brownian motion and the definition of the operator $K_H$,
we deduce
\beqnn
\tilde{B}_t^H\ar=\ar\int_0^tu_sds+\int_0^tK_H(t,s)dW_s\cr
\ar=\ar\int_0^tK_H(t,s)\left[(K_H^{-1}\int_0^\cdot u_rdr)(s)ds+dW_s\right]\cr
\ar=:\ar\int_0^tK_H(t,s)d\tilde{W}_s.
\eeqnn
Now, let
\beqnn
R_T={\rm exp}\left[-\int_0^T\left(K_H^{-1}\int_0^\cdot u_rdr\right)(s)dW_s-\frac{1}{2}\int_0^T\left(K_H^{-1}\int_0^\cdot u_rdr\right)^2(s)ds\right].
\eeqnn
Next we want to show $(\tilde{B}_t^H)_{0\leq t\leq T}$ is an $\mathscr{F}_t^{B^H}$-fractional Brownian motion with Hurst parameter $H$ under the
new probability $R_TP$.
Due to \cite[Theorem 2]{Nualart}, it suffices to show that $\mathbb{E}R_T=1$.
Since $\int_0^\cdot u_rdr$ is absolutely continuous, then
\beqnn
(K_H^{-1}\int_0^\cdot u_rdr)(s)=s^{H-\frac{1}{2}}I_{0+}^{\frac{1}{2}-H}s^{\frac{1}{2}-H}u_s.
\eeqnn
Hence, we have
\beqnn
\left|(K_H^{-1}\int_0^\cdot u_rdr)(s)\right|\ar=\ar\left|\frac{1}{\Gamma(\frac{1}{2}-H)}s^{H-\frac{1}{2}}\int_0^sr^{\frac{1}{2}-H}u_r(s-r)^{-H-\frac{1}{2}}dr\right|\cr
\ar\leq\ar\frac{1}{\Gamma(\frac{1}{2}-H)}s^{H-\frac{1}{2}}\int_0^s|\eta_r|r^{\frac{1}{2}-H}(s-r)^{-H-\frac{1}{2}}dr\cr
\ar\leq\ar\frac{1}{\Gamma(\frac{1}{2}-H)}\frac{|x-y|}{(2K)^{-1}(1-e^{-2KT})}s^{H-\frac{1}{2}}\int_0^sr^{\frac{1}{2}-H}(s-r)^{-H-\frac{1}{2}}dr\cr
\ar=\ar\frac{B(\frac{3}{2}-H,\frac{1}{2}-H)}{\Gamma(\frac{1}{2}-H)} \frac{|x-y|}{(2K)^{-1}(1-e^{-2KT})} s^{\frac{1}{2}-H}.
\eeqnn
As a consequence, we get
\beqlb\label{3.4}
\mathbb{E}{\rm exp}\left[\frac{1}{2}\int_0^T\left(K_H^{-1}\int_0^\cdot u_rdr\right)^2(s)ds\right]
\leq {\rm exp}[C(T,K,H)|x-y|^2],
\eeqlb
where $C(T,K,H)=\left(\frac{B(\frac{3}{2}-H,\frac{1}{2}-H)}{\Gamma(\frac{1}{2}-H)}\right)^2\frac{T^{2-2H}}{K^{-2}(1-e^{-2KT})^2(1-H)}$.
Using the Novikov criterion, we have $\mathbb{E}R_T=1$.\\
Step 3. From step 2, we can rewrite \eqref{3.1} in the following form
\beqnn
dY_t=b(t,Y_t)dt+d\tilde{B}_t^H,\ \ Y_0=y,
\eeqnn
where $(\tilde{B}_t^H)_{0\leq t\leq T}$ is an $\mathscr{F}_t^{B^H}$-fractional Brownian motion with Hurst parameter $H$ under the
new probability $R_TP$.
By the uniqueness of the solution and $X_T=Y_T, a.s.$, we have
\beqlb\label{3.5}
P_Tf(y)=\mathbb{E}f(X_T^y)=\mathbb{E}R_Tf(Y_T^y)=\mathbb{E}R_Tf(X_T^x).
\eeqlb
Applying the H\"{o}lder inequality to \eqref{3.5}, we obtain
\beqlb\label{3.6}
(P_Tf(y))^p\leq P_Tf^p(x)\cdot(\mathbb{E}R_T^\frac{p}{p-1})^{p-1}.
\eeqlb
Now we will estimate moments of $R_T$.\\
Denote $\alpha=\frac{p}{p-1}$ and $M_T=-\int_0^T\left(K_H^{-1}\int_0^\cdot u_rdr\right)(s)dW_s$. Since $(R_t)_{0\leq t\leq T}$ is a $\mathbb{P}$ martingale, by \eqref{3.4} we have
\beqlb\label{3.7}
\mathbb{E}R_T^\alpha\ar=\ar\mathbb{E}{\rm exp}[\alpha M_T-\frac{1}{2}\alpha\langle M\rangle_T]\cr
\ar=\ar\mathbb{E}{\rm exp}[\alpha M_T-\frac{1}{2}\alpha^2\langle M\rangle_T+\frac{1}{2}\alpha(\alpha-1)\langle M\rangle_T]\cr
\ar\leq\ar{\rm exp}[\alpha(\alpha-1)C(T,K,H)|x-y|^2].
\eeqlb
Substituting \eqref{3.7} into \eqref{3.6}, we get the desired result.

\bremark\label{r3.1}
In the proof of Theorem \ref{t3.1}, the choice of $u_t$ is not unique. For instance, we can take another as follows
\beqnn
u_t=\eta_t\cdot\frac{X_t-Y_t}{|X_t-Y_t|}I_{\{t<\tau\}},
\eta_t=\frac{1}{\int_0^Te^{-Kt}dt}\cdot|x-y|.
\eeqnn
Correspondingly, the result of Theorem \ref{t3.1} is of the following form
\beqnn
(P_Tf(y))^p\leq P_Tf^p(x){\rm exp}[\frac{p}{p-1}\tilde{C}(T,K,H)|x-y|^2],
\eeqnn
where
$\tilde{C}(T,K,H)=\left(\frac{B(\frac{3}{2}-H,\frac{1}{2}-H)}{\Gamma(\frac{1}{2}-H)}\right)^2\frac{T^{2-2H}}{4K^{-2}(1-e^{-KT})^2(1-H)}$.
\eremark

As applications of Theorem \ref{t3.1}, we prove the following results on strong Feller property for $P_T$ and log-Harnack inequality.
\bproposition\label{p3.1} Assume {\rm(H1)}. Then $P_T$ is strong feller and the
following estimate holds
\beqnn|P_Tf(x)-P_Tf(y)|\leq ||f||_\infty [2C(T,K,H)]^\frac{1}{2}|x-y|{\rm exp}[C(T,K,H)|x-y|^2], \eeqnn
for every $T>0, x,y\in\mathbb{R}$ and $f\in\mathscr{B}_b(\mathbb{R} )$.
\eproposition

\emph{Proof.} It follows from the proof of Theorem \ref{t3.1} that, for each $f\in\mathscr{B}_b(\mathbb{R} )$,
\beqlb\label{3.8}
|P_Tf(x)-P_Tf(y)|=|\mathbb{E}f(X_T^x)-\mathbb{E}R_Tf(X_T^x)|\leq||f||_\infty\mathbb{E}|1-R_T|.
\eeqlb
Next we will estimate the term $\mathbb{E}|1-R_T|$.\\
Firstly, we have
\beqlb\label{3.9}
(\mathbb{E}|1-R_T|)^2\leq \mathbb{E}|1-R_T|^2=\mathbb{E}R_T^2-1.
\eeqlb
Taking $\alpha=2$ in \eqref{3.7}, we have
\beqlb\label{3.10}
\mathbb{E}R_T^2\leq{\rm exp}[2C(T,K,H)|x-y|^2].
\eeqlb
Combining \eqref{3.9} with \eqref{3.10}, we get
\beqlb\label{3.11}
(\mathbb{E}|1-R_T|)^2\leq 2C(T,K,H)|x-y|^2{\rm exp}[2C(T,K,H)|x-y|^2],
\eeqlb
where we use the elementary inequality $e^x-1\leq xe^x,\ \forall x\geq0$.\\
Substituting \eqref{3.11} into \eqref{3.8}, we can deduce the desired result.

\bcorollary\label{c3.1}
Let {\rm(H1)} hold, then
\beqnn
P_T(\log f)(x)\leq \log P_Tf(y)+C(T,K,H)|x-y|^2,
\eeqnn
$\forall x,y\in\mathbb{R}, t>0, f\geq 1, f\in \mathscr{B}_b(\mathbb{R})$.
\ecorollary
That is, log-Harnack inequality holds.\\
In fact, since $\mathbb{R}$ is a length space, then, by \cite[Proposition 2.2]{Wang1}, we know the result holds.\\

To state further application of Theorem \ref{t3.1}, let us introduce another assumption and some notations.\\
(H2): let $\mu$ be a probability measure on $\mathbb{R}$ such that for some $\tilde{K}>0$,
\beqnn
\mu(P_Tf)\leq \tilde{K}\mu(f), \forall f\in \mathscr{B}_b^+(\mathbb{R}).
\eeqnn
 Note that if $\mu$ is $P_T$-invariant, then (H2) holds.
\bremark\label{r3.2}
 The measures $\mu$ satisfying {\rm(H2)} always exist. For instance,
\beqnn
\mu=\sum\limits_{n=1}^\infty\frac{1}{2^n}P_T^n(x,\cdot),\ \forall x\in \mathbb{R},
\eeqnn
where $(P_T^n(x,\cdot))_{n\geq 1}$ is defined recursively as follows
\beqnn P_T(x,A):=P_TI_A(x),\ \
P_T^n(x,A):=\int_\mathbb{R}P_T^{n-1}(x,dy)P_T(y,A), \ n\geq 2.
\eeqnn
\eremark
Let $\mathscr{C}(\mu,\nu)$ denote the set of all couplings of $\mu$ and $\nu$, where $\mu$ and $\nu$ are two given probability on $\mathbb{R}$,
and $W_2(\mu,\nu)$ be the $L^2$-Wasserstein distance between them with respect to the Euclidian distance, i.e.
\beqnn
W_2^2(\mu,\nu)=\inf\limits_{\pi\in\mathscr{C}(\mu,\nu)}\int_\mathbb{R}\int_\mathbb{R}|x-y|^2\pi(dx,dy).
\eeqnn

\bcorollary\label{c3.2}
Assume that {\rm(H1)} holds and $\mu$ satisfies {\rm(H2)} $(\tilde{K}=1)$. Then the following entropy-cost inequality holds for each $T>0$ and $f\in
\mathscr{B}_b^+(\mathbb{R})$ with $\mu(f)=1$,
\beqnn
\mu(P_T^*f\log P_T^*f)\leq C(T,K,H)W_2^2(\mu,f\mu),
\eeqnn
where $P_T^*$ is the adjoint operator of $P_T$ in $L^2(\mu)$.
\ecorollary
\emph{Proof.} By Corollary \ref{c3.1} for $P_T^*f$,  we have
\beqlb\label{3.12}
P_T(\log P_T^*f)(x)\leq \log P_T(P_T^*f)(y)+C(T,K,H)|x-y|^2, \ \forall x,y\in\mathbb{R}.
\eeqlb
Integrating both sides of \eqref{3.12} with respect to $\pi\in\mathscr{C}(\mu,f\mu)$, we get
\beqnn
 \mu(P_T^*f\log P_T^*f)\leq \mu( \log P_T(P_T^*f))+C(T,K,H)\int_\mathbb{R}\int_\mathbb{R}|x-y|^2\pi(dx,dy).
\eeqnn
Note that, the Jensen inequality and the hypotheses imply
\beqnn
\mu( \log P_T(P_T^*f))\leq\log\mu(P_T(P_T^*f))\leq\log\mu(P_T^*f)=\log\mu(fP_T1)=\log\mu(f)=0.
\eeqnn
So, we get
\beqnn
 \mu(P_T^*f\log P_T^*f)\leq  C(T,K,H)\inf\limits_{\pi\in\mathscr{C}(\mu,f\mu)}\int_\mathbb{R}\int_\mathbb{R}|x-y|^2\pi(dx,dy).
\eeqnn
The proof is complete.

\section{Derivative formula}

\setcounter{equation}{0}
In this part, we begin with the following hypothesis (H3):
\begin{itemize}
 \item[(i)]
 $\partial_2b(t,x)\leq \overline{K}$,  $\forall x\in\mathbb{R}, t\in
 [0,T]$, where $\overline{K}>0$ is a constant, where $\partial_2b(t,x)$ denotes the derivative for the second variable; 
 \par
 \item[(ii)]The mapping $t\mapsto b(t,0)$ is bounded on $[0,T]$.
\end{itemize}

The aim is to establish a Bismut type derivative formula for $P_T$ which will imply the Harnack inequality. For $f\in\mathscr{B}_b(\mathbb{R}), x, y\in\mathbb{R}, T>0,$ we will consider
$$D_yP_Tf(x):=\lim\limits_{\epsilon\rightarrow 0}\frac{P_Tf(x+\epsilon y)-P_Tf(x)}{\epsilon}.$$

\btheorem\label{t4.1} (Derivative formula)
Assume {\rm(H3)}. Then, for each $T>0, f\in\mathscr{B}_b(\mathbb{R}), x, y\in\mathbb{R}, D_yP_Tf(x)$ exists and satisfies
$$D_yP_Tf(x)=\mathbb{E}f(X_T^x)N_T,$$
where $N_T=\frac{1}{\Gamma(\frac{1}{2}-H)T}\int_0^Ts^{H-\frac{1}{2}}\left[\int_0^s\frac{r^{\frac{1}{2}-H}}{(s-r)^{\frac{1}{2}+H}}(1+\partial_2b(r,X_r)(T-r))dr\right]ydW_s.$
\etheorem

\emph{Proof.}
As above, $X_t^x$ solves the equation \eqref{2.1}. For any $\epsilon>0$ and $y\in\mathbb{R}$, we introduce the following equation
\beqlb\label{4.1}
dX_t^\epsilon=b(t,X_t)dt+dB_t^H-\frac{\epsilon}{T}y,\  X_0^\epsilon=x+\epsilon y.
\eeqlb
By {\rm (H3)}, we easily know that the above equation has a unique solution.
Combining \eqref{2.1} with \eqref{4.1}, we deduce that $X_t^\epsilon-X_t=\frac{T-t}{T}\epsilon y,\ \forall t\in[0,T]$, in particular, 
$X_T^\epsilon=X_T$.
Let $\eta_t=b(t,X_t)-b(t,X_t^\epsilon)-\frac{\epsilon}{T}y,\ \forall t\in[0,T]$, then we can rewrite \eqref{4.1} in the form:
\beqnn
dX_t^\epsilon=b(t,X_t^\epsilon)dt+d\overline{B}_t^H,
\eeqnn
where $\overline{B}_t^H=B_t^H+\int_0^t\eta_sds.$
Note that
$$|\eta_t|\leq \overline{K}|X_t-X_t^\epsilon|+\frac{\epsilon y}{T}=\overline{K}\frac{T-t}{T}\epsilon y+\frac{\epsilon y}{T},$$
so, we have $\int_0^T\eta_t^2dt\leq\infty$, and moreover, $\int_0^\cdot \eta_rdr\in I_{0+}^{H+\frac{1}{2}}(L^2([0,T]))$.
Due to the integral representation of fractional Brownian motion and the definition of the operator $K_H$, we get
\beqnn
\overline{B}_t^H=\int_0^tK_H(t,s)d\overline{W}_s,
\eeqnn
where $\overline{W}_t=W_t+\int_0^t(K_H^{-1}\int_0^\cdot \eta_rdr)(s)ds$.
Now, let
\beqnn
R_\epsilon={\rm exp}\left[-\int_0^T\left(K_H^{-1}\int_0^\cdot \eta_rdr\right)(s)dW_s-\frac{1}{2}\int_0^T\left(K_H^{-1}\int_0^\cdot \eta_rdr\right)^2(s)ds\right]
\eeqnn
Now we will prove that $(\overline{B}_t^H)_{0\leq t\leq T}$ is an $\mathscr{F}_t^{B^H}$-fractional Brownian motion with Hurst parameter $H$ under the
new probability $R_\epsilon P$, according to \cite[Theorem 2]{Nualart}, it only needs to show $\mathbb{E}R_\epsilon=1$. 
Similar to step 2 of theorem \ref{t3.1}, we get
\beqlb\label{4.2}
\left|(K_H^{-1}\int_0^\cdot\eta_rdr)(s)\right|
\leq\frac{B(\frac{3}{2}-H,\frac{1}{2}-H)}{\Gamma(\frac{1}{2}-H)}\epsilon y(\overline{K}+\frac{1}{T})s^{\frac{1}{2}-H}.
\eeqlb
Hence, it follows that
\beqnn
\mathbb{E}{\rm exp}\left[\frac{1}{2}\int_0^T\left(K_H^{-1}\int_0^\cdot \eta_rdr\right)^2(s)ds\right]
<\infty.
\eeqnn
By the Novikov criterion, $\mathbb{E}R_T=1$ holds.\\
Hence, in view of the uniqueness of the solution and $X_T^\epsilon=X_T$, we have
\beqnn
P_Tf(x+\varepsilon y)=\mathbb{E}R_\epsilon f(X_t^x).
\eeqnn
With the help of the dominated convergence theorem due to \eqref{4.2}, we deduce that
\beqnn
D_yP_Tf(x):\ar=\ar\lim\limits_{\epsilon\rightarrow 0}\frac{P_Tf(x+\epsilon y)-P_Tf(x)}{\epsilon}\cr
\ar=\ar\lim\limits_{\epsilon\rightarrow 0}\left[\mathbb{E}f(X_T^x)\frac{R_\epsilon-1}{\epsilon}\right]\cr
\ar=\ar\mathbb{E}\left[f(X_T^x)\lim\limits_{\epsilon\rightarrow 0}\frac{R_\epsilon-1}{\epsilon}\right].\cr
\eeqnn
Let $\widetilde{M}_T=:-\int_0^T\left(K_H^{-1}\int_0^\cdot \eta_rdr\right)(s)dW_s.$
Thanks to \eqref{4.2}, we get
\beqnn
\langle \widetilde{M}\rangle_T=\int_0^T|(K_H^{-1}\int_0^\cdot\eta_rdr)(s)|^2ds\leq C\epsilon^2,
\eeqnn
where $C$ is a positive constant.
Therefore, we deduce
\beqnn
\lim\limits_{\epsilon\rightarrow 0}\frac{R_\epsilon-1}{\epsilon}
\ar=\ar\lim\limits_{\epsilon\rightarrow 0}\frac{{\rm exp}[\widetilde{M}_T-\frac{1}{2}\langle \widetilde{M}\rangle_T]-1}{\epsilon}\cr
\ar=\ar\lim\limits_{\epsilon\rightarrow 0}\frac{\widetilde{M}_T-\frac{1}{2}\langle\widetilde{M}\rangle_T}{\epsilon}\cr
\ar=\ar\lim\limits_{\epsilon\rightarrow 0}\frac{\widetilde{M}_T}{\epsilon}.\cr
\eeqnn
Note that
\beqnn
\widetilde{M}_T\ar=\ar-\int_0^T\left(K_H^{-1}\int_0^\cdot \eta_rdr\right)(s)dW_s\cr
\ar=\ar-\frac{1}{\Gamma{(\frac{1}{2}-H})}\int_0^Ts^{H-\frac{1}{2}}\int_0^s\frac{r^{\frac{1}{2}-H}}{(s-r)^{H+\frac{1}{2}}}\eta_rdrdW_s\cr
\ar=\ar\frac{1}{\Gamma{(\frac{1}{2}-H})}\int_0^Ts^{H-\frac{1}{2}}\int_0^s\frac{r^{\frac{1}{2}-H}}{(s-r)^{H+\frac{1}{2}}}[b(r,X_r^\epsilon)-b(r,X_r]drdW_s\cr
\ar+\ar\epsilon\frac{1}{\Gamma{(\frac{1}{2}-H})}\int_0^Ts^{H-\frac{1}{2}}\int_0^s\frac{r^{\frac{1}{2}-H}}{(s-r)^{H+\frac{1}{2}}}\frac{y}{T}drdW_s,\cr
\eeqnn
therefore by {\rm (H3)}, we conclude that
\beqnn
\lim\limits_{\epsilon\rightarrow 0}\frac{R_\epsilon-1}{\epsilon}
=\frac{1}{\Gamma(\frac{1}{2}-H)T}\int_0^Ts^{H-\frac{1}{2}}\left[\int_0^s\frac{r^{\frac{1}{2}-H}}{(s-r)^{\frac{1}{2}+H}}(1+\partial_2b(r,X_r)(T-r))dr\right]ydW_s.
\eeqnn 
The proof is complete.

\bremark\label{r4.1}
If H=$\frac{1}{2}$, i.e. $B^H$ is a Brownian motion, then the corresponding derivative formula is of the following type
$$D_yP_Tf(x)=\frac{1}{T}\int_0^T[1+(T-s)\partial_2b(s,X_s)]ydW_s.$$
\eremark

\bremark\label{r4.2}
Since we deal with one dimensional case, the derivative formula of theorem \eqref{4.1} is equivalent to $P_Tf(\cdot)$ is derivative. The method we adopt is also valid for $n$-dimensional case.
\eremark

As an application of the  derivative formula derived above, we have the following result.

\bcorollary\label{c4.1}
If {\rm(H3)} holds, then for any nonnegative $f\in\mathscr{B}_b(\mathbb{R})$ and $t>0, x,y\in\mathbb{R}$,
\beqnn
(P_Tf(y))^p\leq P_Tf^p(x){\rm exp}[\frac{p}{p-1}C(T,\overline{K},H)|y-x|^2],
\eeqnn
where
$C(T,\overline{K},H)=\left(\frac{B(\frac{3}{2}-H,\frac{1}{2}-H)}{\Gamma(\frac{1}{2}-H)}\right)^2\frac{(1+\overline{K}T)^2}{T^{2H}2(1-H)}$.
\ecorollary

\emph{Proof.} By \eqref{4.1} and the Young inequality \cite{Arnaudon2}, we have. for all $\delta>0$,
\beqlb\label{4.3}
|D_yP_Tf(x)|\leq\delta[P_T(f\log f)(x)-(P_Tf)(x)(\log P_Tf)(x)]+P_Tf(x)[\delta\log\mathbb{E}e^{\frac{1}{\delta}N_T}].
\eeqlb
 Now let $\beta(s)=1+s(p-1), \gamma(s)=x+s(y-x), s\in[0,1]$, we have
\beqnn
\ar\frac{d}{ds}\ar\log (P_Tf^{\beta(s)})^\frac{p}{\beta(s)}(\gamma(s))\cr
\ar=\ar\frac{p(p-1)}{\beta^2(s)}\frac{P_T(f^{\beta(s)}\log f^{\beta(s)})-(P_Tf^{\beta(s)})\log P_Tf^{\beta(s)}}{P_Tf^{\beta(s)}}(\gamma(s))
+\frac{p}{\beta(s)}\frac{D_{y-x}P_Tf^{\beta(s)}}{P_Tf^{\beta(s)}}(\gamma(s))\cr
\ar\geq\ar\frac{p}{\beta(s)P_Tf^{\beta(s)}(\gamma(s))} \{\frac{p-1}{\beta(s)}[P_T(f^{\beta(s)}\log f^{\beta(s)})(\gamma(s))\cr
\ar-\ar(P_Tf^{\beta(s)})\log P_Tf^{\beta(s)}(\gamma(s))]-|D_{y-x}P_Tf^{\beta(s)}|(\gamma(s)\}\cr
\ar\geq\ar-\frac{p(p-1)}{\beta^2(s)}\log \mathbb{E}e^{\frac{1}{\delta}N_T},\cr
\eeqnn
where we use \eqref{4.3} and choose $\delta=\frac{p-1}{\beta(s)}$ for the last inequality, note that $N_T$ is corresponding to the direction $y-x$.\\
Next we are to estimate $\mathbb{E}e^{\frac{1}{\delta}N_T}$. Since
 $\mathbb{E}e^{\frac{1}{\delta}N_T}\leq\left(\mathbb{E}e^{\frac{2}{\delta^2}\langle N_T\rangle}\right)^\frac{1}{2}$, we turn to
the term $\langle N_T\rangle$.
\beqnn
\langle N_T\rangle\ar=\ar
\frac{1}{(\Gamma(\frac{1}{2}-H)T)^2}\int_0^Ts^{2H-1}\left[\int_0^s\frac{r^{\frac{1}{2}-H}}{(s-r)^{\frac{1}{2}+H}}(1+\partial_2b(r,X_r)(T-r))dr\right]^2(y-x)^2ds\cr
\ar\leq\ar\left(\frac{1+\overline{K}T}{\Gamma(\frac{1}{2}-H)T}\right)^2\int_0^Ts^{2H-1}\left[\int_0^s\frac{r^{\frac{1}{2}-H}}{(s-r)^{\frac{1}{2}+H}}dr\right]^2(y-x)^2ds\cr
\ar=\ar\left(\frac{B(\frac{3}{2}-H,\frac{1}{2}-H)}{\Gamma(\frac{1}{2}-H)}\right)^2\frac{(1+\overline{K}T)^2}{T^{2H}2(1-H)}(y-x)^2\cr
\ar=:\ar C(T,\overline{K},H)(y-x)^2.\cr
\eeqnn
Therefore, we deduce that
\beqnn
\ar\frac{d}{ds}\ar\log (P_Tf^{\beta(s)})^\frac{p}{\beta(s)}(\gamma(s))\cr
\ar\geq\ar-\frac{p}{p-1}C(T,\overline{K},H)(y-x)^2.
\eeqnn
Integrating on the interval $[0,1]$ with respect to $s$, we get the desired result.

\bremark\label{r4.4}
To our knowledge, for the results of theorem \ref{t3.1} and corollary \ref{c4.1}, we can not decide which is better.
\eremark

\end{document}